%
%
%
%
%
%
\magnification=\magstephalf      
%
%
\vsize=7.5truein                 
\hsize=5.2truein                 
\newskip\stdskip                 
\stdskip=6pt plus3pt minus3pt    
\medskipamount=\stdskip          
\parindent=0pt                   
\parskip=\stdskip                
\abovedisplayskip=\stdskip       
\belowdisplayskip=\stdskip       
\mathsurround=0.75pt             
\overfullrule=0pt                
%
%
\def\ppar{\par\goodbreak\vskip 8pt plus 4pt minus 4pt}     
%
%
\def\stdspace{\hskip 0.75em plus 0.15em\ignorespaces}
%
%
%
%
%
%
%
\def\hexnumber#1{\ifcase#1 0\or 1\or 2\or 3\or 4\or 5\or 6\or 7\or 8\or
 9\or A\or B\or C\or D\or E\or F\fi}
%
%
\font\thirtnmsa=msam10 scaled 1315    
\font\tenmsa=msam10          \font\ninemsa=msam9
\font\sevenmsa=msam7         \font\sixmsa=msam6
\font\fivemsa=msam5
%
%
\newfam\msafam                  \textfont\msafam=\tenmsa
\scriptfont\msafam=\sevenmsa    \scriptscriptfont\msafam=\fivemsa
\edef\hexa{\hexnumber\msafam}        
\def\msa{\fam\msafam\tenmsa}         
%
%
\font\thirtnmsb=msbm10 scaled 1315   
\font\tenmsb=msbm10      \font\ninemsb=msbm9
\font\sevenmsb=msbm7     \font\sixmsb=msbm6
\font\fivemsb=msbm5
%
\newfam\msbfam                   \textfont\msbfam=\tenmsb       
\scriptfont\msbfam=\sevenmsb     \scriptscriptfont\msbfam=\fivemsb
\edef\hexb{\hexnumber\msbfam}    
\def\msb{\fam\msbfam\tenmsb}     
%
%
\font\thirtneufm=eufm10 scaled 1315   
\font\teneufm=eufm10                 \font\nineeufm=eufm9
\font\seveneufm=eufm7                \font\sixeufm=eufm6
\font\fiveeufm=eufm5
%
\newfam\eufmfam                    \textfont\eufmfam=\teneufm
\scriptfont\eufmfam=\seveneufm     \scriptscriptfont\eufmfam=\fiveeufm
\edef\hexf{\hexnumber\eufmfam}      
\def\frak{\fam\eufmfam\teneufm}     
%
%
%
\font\thirtnrm=cmr10 scaled 1315    
\font\ninerm=cmr9                   \font\sixrm=cmr6   
%
\font\thirtni=cmmi10 scaled 1315    
\font\ninei=cmmi9                   \font\sixi=cmmi6  
%
\font\thirtnsy=cmsy10 scaled 1315   
\font\ninesy=cmsy9                  \font\sixsy=cmsy6  
%
\font\thirtnbf=cmbx10 scaled 1315   
\font\ninebf=cmbx9                  \font\sixbf=cmbx6  
%
%
\font\thirtnex=cmex10 scaled 1315   
\font\nineex=cmex9                  
%
%
\font\thirtnit=cmti10 scaled 1315  
\font\nineit=cmti9                  
%
\font\thirtnsl=cmsl10 scaled 1315  
\font\ninesl=cmsl9                  
%
\font\thirtntt=cmtt10 scaled 1315  
\font\ninett=cmtt9                  
%
%
%
%
\def\small{%
%
%
\textfont0=\ninerm \scriptfont0=\sixrm \scriptscriptfont0=\fiverm
\def\rm{\fam0\ninerm}
%
%
\textfont1=\ninei \scriptfont1=\sixi \scriptscriptfont1=\fivei
%
%
\textfont2=\ninesy \scriptfont2=\sixsy \scriptscriptfont2=\fivesy
%
%
\textfont3=\nineex \scriptfont3=\nineex \scriptscriptfont3=\nineex
%
%
\textfont\bffam=\ninebf \scriptfont\bffam=\sixbf
\scriptscriptfont\bffam=\fivebf \def\bf{\fam\bffam\ninebf}%
%
%
\textfont\itfam=\nineit \def\it{\fam\itfam\nineit}%
\textfont\slfam=\ninesl \def\sl{\fam\slfam\ninesl}%
\textfont\ttfam=\ninett \def\tt{\fam\ttfam\ninett}%
%
%
%
\textfont\msafam=\ninemsa \scriptfont\msafam=\sixmsa
\scriptscriptfont\msafam=\fivemsa \def\msa{\fam\msafam\ninemsa}%
%
%
\textfont\msbfam=\ninemsb \scriptfont\msbfam=\sixmsb
\scriptscriptfont\msbfam=\fivemsb \def\msb{\fam\msbfam\ninemsb}%
%
%
\textfont\eufmfam=\nineeufm  \scriptfont\eufmfam=\sixeufm
\scriptscriptfont\eufmfam=\fiveeufm \def\frak{\fam\eufmfam\nineeufm}%
%
%
%
\normalbaselineskip=11pt%
\setbox\strutbox=\hbox{\vrule height8pt depth3pt width0pt}%
%
%
\normalbaselines\rm
%
%
\stdskip=4pt plus2pt minus2pt    
\medskipamount=\stdskip          
\parskip=\stdskip                
\abovedisplayskip=\stdskip       
\belowdisplayskip=\stdskip       
\def\ppar{\par\goodbreak\vskip 6pt plus 3pt minus 3pt}%
%
%
\def\section##1{\global\advance\sectionnumber by 1
\vskip-\lastskip\penalty-800\vskip 20pt plus10pt minus5pt 
\egroup{\bf\number\sectionnumber\quad##1}\bgroup\small         
\vskip 6pt plus3pt minus3pt
\nobreak\resultnumber=1}
}    
%
\def\beginsmall{\bgroup\small}
\let\endsmall\egroup
%
%
%
%
\def\large{%
\textfont0=\thirtnrm \scriptfont0=\ninerm \scriptscriptfont0=\sevenrm
\def\rm{\fam0\thirtnrm}%
\textfont1=\thirtni \scriptfont1=\ninei \scriptscriptfont1=\seveni
\textfont2=\thirtnsy \scriptfont2=\ninesy \scriptscriptfont2=\sevensy
\textfont3=\thirtnex \scriptfont3=\thirtnex \scriptscriptfont3=\thirtnex
\textfont\bffam=\thirtnbf \scriptfont\bffam=\ninebf
\scriptscriptfont\bffam=\sevenbf \def\bf{\fam\bffam\thirtnbf}%
\textfont\itfam=\thirtnit \def\it{\fam\itfam\thirtnit}%
\textfont\slfam=\thirtnsl \def\sl{\fam\slfam\thirtnsl}%
\textfont\ttfam=\thirtntt \def\tt{\fam\ttfam\thirtntt}%
\textfont\msafam=\thirtnmsa \scriptfont\msafam=\ninemsa
\scriptscriptfont\msafam=\sevenmsa \def\msa{\fam\msafam\thirtnmsa}%
\textfont\msbfam=\thirtnmsb \scriptfont\msbfam=\ninemsb
\scriptscriptfont\msbfam=\sevenmsb \def\msb{\fam\msbfam\thirtnmsb}%
\textfont\eufmfam=\thirtneufm  \scriptfont\eufmfam=\nineeufm
\scriptscriptfont\eufmfam=\seveneufm \def\frak{\fam\eufmfam\teneufm}%
\normalbaselineskip=16pt%
\setbox\strutbox=\hbox{\vrule height11.5pt depth4.5pt width0pt}%
\normalbaselines\rm}%
\let\Large\large   
%

%

%
\mathchardef\plussquare="0\hexa01
\mathchardef\nge="3\hexb0B
\mathchardef\maltesecross="0\hexa7A
\mathchardef\del="0\hexf01
%
%
%
%
\font\sc=cmcsc10
%
%
%
%
\def\sqr#1#2{{\vcenter{\vbox{\hrule  height.#2truept
	\hbox{\vrule width.#2truept height#1truept 
	\kern#1truept \vrule width.#2truept}
	\hrule height.#2truept}}}}
%
%
%
%
%
\newcount\sectionnumber            
\newcount\resultnumber             
\sectionnumber=0\resultnumber=1    
%
%
%
\def\section#1{\global\advance\sectionnumber by 1
\xdef\nextkey{\number\sectionnumber}
\vskip-\lastskip\penalty-800\vskip 20pt plus10pt minus5pt 
{\large\bf\number\sectionnumber\quad#1}         
\vskip 8pt plus4pt minus4pt
\nobreak\resultnumber=1}      
%
%
%
%
%
         
%
%
%
%

%
\def\proc#1{\xdef\nextkey{\number\sectionnumber.\number\resultnumber}%
\vskip-\lastskip\ppar\bf%
\noindent#1\ \number\sectionnumber.\number\resultnumber
\stdspace\sl\global\advance\resultnumber by 1\ignorespaces}
\def\endproc{\rm\ppar} 
%
%
%
%
%
%
%
%
%
%
\def\proclaim#1{\vskip-\lastskip\ppar\bf%
\noindent#1\stdspace\sl\ignorespaces} 
\let\endproclaim\endproc
%
%
%
%

%
%
%
%
%
%
\def\label{\xdef\nextkey{\number\sectionnumber.\number\resultnumber}%
\number\sectionnumber.\number\resultnumber
\global\advance\resultnumber by 1}
%
%
%
%
%
%
%
%
%
%
%
%
%
%
%
%
\newcount\refnumber              
\refnumber=1                     
\long\def\reflist#1\endreflist{%
\long\def\thereflist{#1}{\def\refkey##1##2\par{\xdef##1{\number\refnumber}%
\global\advance\refnumber by 1}%
\def\key##1##2\par{\expandafter\xdef%
\csname##1\endcsname{\number\refnumber}%
\global\advance\refnumber by 1}#1\par}}
\long\def\references{%
\penalty-800\vskip-\lastskip\vskip 15pt plus10pt minus5pt 
{\large\bf References}\ppar 
{\leftskip=25pt\frenchspacing    
\small\parskip=3pt plus2pt       
\def\refkey##1##2\par{\noindent  
\llap{[##1]\stdspace}\ignorespaces##2\par}         
\def\key##1##2\par{\noindent  
\llap{[\ref{##1}]\stdspace}\ignorespaces##2\par}  
\def\,{\thinspace}\thereflist\par}}
%
%
%
\newcount\footnotenumber         
\footnotenumber=1                
\def\fnote#1{\xdef\nextkey{\number\footnotenumber}%
{\small\ifnum\footnotenumber>9\parindent=14pt%
\else\parindent=10pt\fi\footnote{$^{\number\footnotenumber}$}%
{\hglue-5pt#1}\global\advance\footnotenumber by 1}}
%
%
%
%
%
%
%
\newcount\figurenumber          
\figurenumber=1                 
\def\caption#1{\xdef\nextkey{\number\figurenumber}%
\cl{\small Figure \number\figurenumber: #1}%
\global\advance\figurenumber by 1}
\def\figurelabel{\xdef\nextkey{\number\figurenumber}%
\cl{\small Figure \number\figurenumber}%
\global\advance\figurenumber by 1}
\long\def\figure#1\endfigure{{\xdef\nextkey{\number\figurenumber}%
\let\captiontext\relax\def\caption##1{\xdef\captiontext{##1}}%
\midinsert\cl{\ignorespaces#1\unskip\unskip\unskip\unskip}\vglue6pt\cl{\small 
Figure \number\figurenumber\ifx\captiontext\relax\else: \captiontext
\fi}\endinsert\global\advance\figurenumber by 1}}
%
%
%
%
%
%
%
\def\nextkey{??}   
%
\def\key#1{\expandafter\xdef\csname #1\endcsname{\nextkey}}
\def\ref#1{\expandafter\ifx\csname #1\endcsname\relax
\immediate\write16{Reference {#1} undefined}??\else
\csname #1\endcsname\fi}
%
%
%
%
%
%
%
\newread\gtinfile
\newwrite\gtreffile
\def\useforwardrefs{
\openin\gtinfile\jobname.ref
\ifeof\gtinfile
\closein\gtinfile
\immediate\write16{No file \jobname.ref}
\else
\closein\gtinfile
\input \jobname.ref
\fi
\immediate\openout\gtreffile \jobname.ref
%
%
\def\key##1{{\def\\{\noexpand}%
\expandafter\xdef\csname ##1\endcsname{\nextkey}%
\immediate\write\gtreffile{\\\expandafter\\\def\\\csname ##1\\\endcsname%
{\nextkey}}}}
%
%
\long\def\reflist##1\endreflist{%
\long\def\thereflist{##1}{\def\refkey####1####2\par{\xdef####1{%
\number\refnumber}{\def\\{\noexpand}\immediate\write\gtreffile
{\\\def\\####1{\number\refnumber}}}\global\advance\refnumber by 1}%
\def\key####1####2\par{\expandafter\xdef%
\csname####1\endcsname{\number\refnumber}%
{\def\\{\noexpand}\immediate\write\gtreffile
{\\\expandafter\\\def\\\csname ####1\\\endcsname{\number\refnumber}}}
\global\advance\refnumber by 1}##1\par}}
\long\def\biblio##1\endbiblio{\reflist##1\endreflist\references}%
%
%
\def\numkey##1{{\def\\{\noexpand}%
\xdef##1{\number\sectionnumber.\number\resultnumber}
\immediate\write\gtreffile{\\\def\\##1%
{\number\sectionnumber.\number\resultnumber}}}}
\def\seckey##1{{\def\\{\noexpand}\xdef##1{\number\sectionnumber}
\immediate\write\gtreffile{\\\def\\##1{\number\sectionnumber}}}}
\def\figkey##1{\xdef##1{\number\figurenumber}%
{\def\\{\noexpand}\immediate\write\gtreffile%
{\\\def\\##1{\number\figurenumber}}}
\number\figurenumber\global\advance\figurenumber by 1}
}   
%
%
%
%
\def\figkey#1{\xdef#1{\number\figurenumber}%
\number\figurenumber\global\advance\figurenumber by 1}
\def\fig#1#2\endfig{%
\midinsert\cl{#2}\vglue6pt\cl{\small Figure #1}\endinsert}
\def\newfig{\number\figurenumber\global\advance\figurenumber by 1}
\def\numkey#1{\xdef#1{\number\sectionnumber.\number\resultnumber}}
\def\seckey#1{\xdef#1{\number\sectionnumber}}
%
%
%
%
%
%
%
%
%
\def\verb{\catcode`\"=\active}       
\def\brev{\catcode`\"=12}            
\brev                                
\verb                                
{\obeyspaces\gdef {\ }}              
{\catcode`\`=\active\gdef`{\relax\lq}}
\def"{%
\begingroup\baselineskip=12pt\def\par{\leavevmode\endgraf}%
\tt\obeylines\obeyspaces\parskip=0pt\parindent=0pt%
\catcode`\$=12\catcode`\&=12\catcode`\^=12\catcode`\#=12%
\catcode`\_=12\catcode`\~=12%
\catcode`\{=12\catcode`\}=12\catcode`\%=12\catcode`\\=12%
\catcode`\`=\active\let"\endgroup}
\brev      
%
%
%
%
%
%
\def\item#1{\par\leavevmode\llap{#1\stdspace}%
\ignorespaces}                             
%
%

%
%
\def\co{\colon\thinspace}    
\def\np{\vfil\eject}         
\def\nl{\hfil\break}         
\def\cl{\centerline}         
\def\gt{{\mathsurround=0pt\it $\cal G\mskip-2mu$eometry \&\ 
$\cal T\!\!$opology}}        
\def\agt{{\mathsurround=0pt\it$\cal A\mskip-.7mu$lgebraic \&\ 
$\cal G\mskip-2mu$eometric $\cal T\!\!$opology}}  
%
%
%

%
%
%
%
%
\def\title#1{\def\thetitle{#1}}

\def\author#1{\edef\previousauthors{\theauthors}
 \ifx\theauthors\relax\def\theauthors{#1}\else
 \def\theauthors{\previousauthors\par#1}\fi}

\let\authors\author        
\def\address#1{\edef\previousaddresses{\theaddress}
 \ifx\theaddress\relax\def\theaddress{#1}\else
 \def\theaddress{\previousaddresses\par\vskip 2pt\par#1}\fi}
\def\secondaddress#1{\edef\previousaddresses{\theaddress}
 \ifx\theaddress\relax\def\theaddress{#1}\else
 \def\theaddress{\previousaddresses\par{\rm and}\par#1}\fi}   

\def\email#1{\edef\previousemails{\theemail}
 \ifx\theemail\relax\def\theemail{#1}\else
 \def\theemail{\previousemails\hskip 0.75em\relax#1}\fi}
\def\secondemail#1{\edef\previousemails{\theemail}
 \ifx\theemail\relax\def\theemail{#1}\else
 \def\theemail{\previousemails\hskip 0.75em{\rm and}\hskip 0.75em
 \relax#1}\fi}
\def\url#1{\edef\previousurls{\theurl}
 \ifx\theurl\relax\def\theurl{#1}\else
 \def\theurl{\previousurls\hskip 0.75em\relax#1}\fi}
\def\secondurl#1{\edef\previousurls{\theurl}
 \ifx\theurl\relax\def\theurl{#1}\else
 \def\theurl{\previousurls\hskip 0.75em{\rm and}\hskip 0.75em
 \relax#1}\fi}
\long\def\abstract#1\endabstract{\long\def\theabstract{#1}}
\def\primaryclass#1{\def\theprimaryclass{#1}}
\def\secondaryclass#1{\def\thesecondaryclass{#1}}
\def\keywords#1{\def\thekeywords{#1}}
%
%
\let\\\par\let\thetitle\relax\let\theshorttitle\relax
\let\theauthors\relax\let\theshortauthors\relax
\let\theaddress\relax\let\theshortaddress\relax
\let\theemail\relax\let\theurl\relax
\let\theabstract\relax\let\theprimaryclass\relax
\let\thesecondaryclass\relax\let\thekeywords\relax
%
%
%
%
\long\def\maketitlepage{    

\vglue 0.2truein   

%
{\parskip=0pt\leftskip 0pt plus 1fil\def\\{\par\smallskip}{\large
\bf\thetitle}\par\medskip}   

\vglue 0.15truein 

%
{\parskip=0pt\leftskip 0pt plus 1fil\def\\{\par}{\sc\theauthors}
\par\medskip}%
 
\vglue 0.1truein 

%
{\small\parskip=0pt
{\leftskip 0pt plus 1fil\def\\{\par}{\sl\theaddress}\par}
\ifx\theemail\relax\else  
\vglue 5pt \def\\{\stdspace{\rm and}\stdspace} 
\cl{Email:\stdspace\tt\theemail}\fi
\ifx\theurl\relax\else    
\vglue 5pt \def\\{\stdspace{\rm and}\stdspace} 
\cl{URL:\stdspace\tt\theurl}\fi\par}

\vglue 7pt 

{\bf Abstract}

\vglue 5pt

\theabstract

\vglue 7pt 

{\bf AMS Classification numbers}\quad Primary:\quad \theprimaryclass\par

Secondary:\quad \thesecondaryclass

\vglue 5pt 

{\bf Keywords:}\quad \thekeywords

\np  

}    
%
%
\long\def\makeshorttitle{    


%
{\parskip=0pt\leftskip 0pt plus 1fil\def\\{\par\smallskip}{\large
\bf\thetitle}\par\medskip}   

\vglue 0.05truein 

%
{\parskip=0pt\leftskip 0pt plus 1fil\def\\{\par}{\sc\theauthors}
\par\medskip}%
 
\vglue 0.03truein 

%
{\small\parskip=0pt
{\leftskip 0pt plus 1fil\def\\{\par}{\sl\ifx\theshortaddress\relax
\theaddress\else\theshortaddress\fi}\par}
\ifx\theemail\relax\else  
\vglue 5pt \def\\{\stdspace{\rm and}\stdspace} 
\cl{Email:\stdspace\tt\theemail}\fi
\ifx\theurl\relax\else    
\vglue 5pt \def\\{\stdspace{\rm and}\stdspace} 
\cl{URL:\stdspace\tt\theurl}\fi\par}

\vglue 10pt 


{\small\leftskip 25pt\rightskip 25pt{\bf Abstract}\stdspace\theabstract

{\bf AMS Classification}\stdspace\theprimaryclass
\ifx\thesecondaryclass\relax\else; \thesecondaryclass\fi\par
{\bf Keywords}\stdspace \thekeywords\par}
\vglue 7pt
}    
\let\maketitle\makeshorttitle        
%
%

\def\volumenumber#1{\def\thevolumenumber{#1}}
\def\volumeyear#1{\def\thevolumeyear{#1}}
\def\pagenumbers#1#2{\def\startpage{#1}\def\finishpage{#2}}
\def\published#1{\def\publishdate{#1}}
\def\received#1{\def\receiveddate{#1}}

\let\reviseddate\relax
\volumenumber{X}
\volumeyear{20XX}
\pagenumbers{1}{XXX}
\published{XX Xxxember 20XX}

\long\def\makeagttitle{   
\agt\hfill      
\hbox to 60truept{\vbox to 0pt{\vglue -14truept{\bf [Logo here]}\vss}\hss}
\break
{\small Volume \thevolumenumber\ (\thevolumeyear)
\startpage--\finishpage\nl
Published: \publishdate}

\vglue .2truein

{\parskip=0pt\leftskip 0pt plus 1fil\def\\{\par\smallskip}{\large
\bf\thetitle}\par\medskip}   
\vglue 0.05truein 

%
{\parskip=0pt\leftskip 0pt plus 1fil\def\\{\par}{\sc\theauthors}
\par\medskip}%
 
\vglue 0.03truein 


{\small\leftskip 25truept\rightskip 25truept{\bf Abstract}\stdspace\theabstract

{\bf AMS Classification}\stdspace\theprimaryclass
\ifx\thesecondaryclass\relax\else; \thesecondaryclass\fi\par
{\bf Keywords}\stdspace \thekeywords\par}\vglue 7truept

}   


\def\Addresses{\bigskip
{\small \parskip 0pt \leftskip 0pt \rightskip 0pt plus 1fil \def\\{\par}
\sl\theaddress\par\medskip \rm Email:\stdspace\tt\theemail\par
\ifx\theurl\relax\else\smallskip \rm URL:\stdspace\tt\theurl\par\fi}}

\def\agtart{
\hoffset 14truemm
\voffset 31truemm
\font\phead=cmsl9 scaled 950
\font\pnum=cmbx10 scaled 913
\font\pfoot=cmsl9 scaled 950
\headline{\vbox to 0pt{\vskip -4.5mm\line{\small\phead\ifnum
\count0=\startpage ISSN numbers are printed here
\hfill {\pnum\folio}\else\ifodd\count0\def\\{ }%
\ifx\theshorttitle\relax\thetitle\else\theshorttitle\fi\hfill{\pnum\folio}
\else\def\\{ and }{\pnum\folio}\hfill\ifx\theshortauthors\relax\theauthors
\else\theshortauthors\fi\fi\fi}\vss}}
\footline{\vbox to 0pt{\vglue 0mm\line{\small\pfoot\ifnum\count0=\startpage
Copyright declaration is printed here\hfill\else
\agt, Volume \thevolumenumber\ (\thevolumeyear)\hfill\fi}\vss}}
\let\maketitle\makeagttitle\let\makeshorttitle\makeagttitle}


\def\ifplaintex{\expandafter\ifx\csname documentclass\endcsname\relax}


\ifplaintex 
\hoffset 14truemm
\voffset 31truemm
\else
\headsep 23pt
\footskip 35pt
\hoffset -4truemm
\voffset 12.5truemm
\fi

\expandafter\ifx\csname beginpicture\endcsname\relax
\expandafter\ifx\csname documentclass\endcsname\relax
\input pictex \else\font\fiverm=cmr5
\input prepictex \input pictex \input postpictex \fi\fi

\def\gt{{\mathsurround=0pt\it $\cal G\mskip-2mu$eometry \&\ 
$\cal T\!\!$opology}}        

\def\gtp{{\mathsurround=0pt\it $\cal G\mskip-2mu$eometry \&\ 
$\cal T\!\!$opology $\cal P\!$ublications}}  


\def\lognumber#1{\def\thelognumber{#1}}
\def\volumenumber#1{\def\thevolumenumber{#1}}
\def\papernumber#1{\def\thepapernumber{#1}}
\def\volumeyear#1{\def\thevolumeyear{#1}}

\def\pagenumbers#1#2{\def\startpage{#1}\def\finishpage{#2}}
\def\published#1{\def\publishdate{#1}}
\def\proposed#1{\def\theproposer{#1}}
\def\seconded#1{\def\theseconders{#1}}
\def\received#1{\def\receiveddate{#1}}

\def\accepted#1{\def\accepteddate{#1}}
\def\asciititle#1{\def\theasciititle{#1}}

\long\def\asciiabstract#1{\long\def\theasciiabstract{#1}}


\let\\\par\let\thelognumber\relax
\let\thevolumenumber\relax\let\thepapernumber\relax
\let\thevolumeyear\relax\let\thesamplenumber\relax\let\startpage\relax
\let\finishpage\relax\let\publishdate\relax\let\receiveddate\relax
\let\reviseddate\relax\let\accepteddate\relax\let\theasciititle\relax
\let\theasciiauthors\relax
\let\theasciiabstract\relax
\let\theasciiemail\relax\let\theshortauthors\relax\let\theshorttitle\relax

\long\def\maketitlep{   

\count0=\startpage

\gt\hfill      
\beginpicture
\setcoordinatesystem units <0.33truein, 0.33truein> point at 2.2 0.9
\setplotsymbol ({$\cal G$})
\plotsymbolspacing=9truept
\circulararc 315 degrees from 0 1 center at 0 0
\setplotsymbol ({$\cal T$})
\circulararc 315 degrees from 1 -1 center at 1 0
\endpicture
%
\break
{\small\ifx\thesamplenumber\relax 
Volume \else Sample
\fi\thevolumenumber\ (\thevolumeyear)
\startpage--\finishpage\nl
Published: \publishdate}
\vglue 0.5truein plus 0.4fil minus 0.1truein

{\parskip=0pt\leftskip 0pt plus 1fil\def\\{\par\smallskip}{\ifplaintex\large
\else\Large\fi\bf\thetitle}\par\medskip}   

\vglue 0pt plus 0.1fil 

{\parskip=0pt\leftskip 0pt plus 1fil\def\\{\par}{\sc\theauthors}
\par\medskip}

\vglue 0pt plus 0.1fil 

{\small\parskip=0pt\let\newline\\
{\leftskip 0pt plus 1fil\def\\{\par}{\sl\theaddress}\par}
\expandafter\ifx\theemail\relax    
\relax\else\vglue 5pt plus 0.02fil minus 2pt\def\\{\stdspace{\rm 
and}\stdspace} 
\cl{Email:\stdspace\tt\theemail}\fi
\ifx\theurl\relax                  
\relax\else\vglue 5pt plus 0.02fil minus 2pt\def\\{\stdspace{\rm 
and}\stdspace}
\cl{URL:\stdspace\tt\theurl}\fi\par}

\vglue 7pt plus 0.3fil minus 3pt

{\bf Abstract}
\vglue 5pt plus 0.1fil minus 2pt

\theabstract

\vglue 7pt plus 0.3fil minus 3pt

{\bf AMS Classification numbers}\quad Primary:\quad \theprimaryclass

Secondary:\quad \thesecondaryclass

\vglue 5pt plus 0.3fil minus 2pt

{\bf Keywords}\quad \thekeywords

\vglue 10pt plus 0.5fil minus 5pt

{\small  Proposed: \theproposer\hfill Received: \receiveddate\nl
Seconded: \theseconders\hfill 
\ifx\reviseddate\relax                         
Accepted: \accepteddate                        
\else
Revised: \reviseddate                          
\fi}
\eject
}       

\let\maketitlepage\maketitlep
\let\maketitle\maketitlepage


\font\phead=cmsl9 scaled 950
\font\lhead=cmsl9 scaled 1050
\font\pnum=cmbx10 scaled 913
\font\lnum=cmbx10 
\font\pfoot=cmsl9 scaled 950
\font\lfoot=cmsl9 scaled 1050
\ifplaintex
\headline{\vbox to 0pt{\vskip -4.5mm\line{\small\phead\ifnum
\count0=\startpage ISSN 1364-0380 (on line)
1465-3060 (printed) \hfill {\pnum\folio}\else\ifodd\count0\def\\{ }%
\ifx\theshorttitle\relax\thetitle\else\theshorttitle\fi\hfill{\pnum\folio}
\else\def\\{ and }{\pnum\folio}\hfill\ifx\theshortauthors\relax\theauthors
\else\theshortauthors\fi\fi\fi}\vss}}
\footline{\vbox to 0pt{\vglue 0mm\line{\small\pfoot\ifnum\count0=\startpage
\copyright\ \gtp\hfill\else
\gt, Volume \thevolumenumber\ (\thevolumeyear)\hfill\fi}\vss
}}
\else
\makeatletter
\def\@oddhead{{\small\lhead\ifnum\count0=\startpage ISSN 1364-0380 (on line)
1465-3060 (printed) \hfill {\lnum\number\count0}\else\ifodd\count0
\def\\{ }\ifx\theshorttitle\relax \thetitle \else\theshorttitle\fi\hfill
{\lnum\number\count0}\else\def\\{ and }{\lnum\number\count0}
\hfill\ifx\theshortauthors\relax 
\theauthors\else\theshortauthors\fi\fi\fi}}\def\@evenhead{\@oddhead}
\def\@oddfoot{\small\lfoot\ifnum\count0=\startpage\copyright\ \gtp\hfill\else
\gt, Volume \thevolumenumber\ (\thevolumeyear)\hfill\fi}
\def\@evenfoot{\@oddfoot}
\makeatother
\fi


\newwrite\gtoutfile
\long\gdef\makeheadfile{  
{\def\\{, }\def\s{ }
\immediate\openout\gtoutfile head.xxx
\immediate\write\gtoutfile{To: math@arxiv.org}
\immediate\write\gtoutfile{Subject: put or rep NNNNN:pppp}
\immediate\write\gtoutfile{--text follows this line--}
\immediate\write\gtoutfile{Proxy-for: \ifx\theasciiauthors\relax
\theauthors\else\theasciiauthors\fi\s<\ifx\theasciiemail\relax\theemail\else\theasciiemail\fi>}
\immediate\write\gtoutfile{\noexpand\\}
\immediate\write\gtoutfile{Authors: \ifx\theasciiauthors\relax
\theauthors\else\theasciiauthors\fi}
{\def\\{ }\immediate\write\gtoutfile{Title: \ifx\theasciititle\relax
\thetitle\else\theasciititle\fi}}
\immediate\write\gtoutfile{Subj-class: GT or SG or MG etc}
\immediate\write\gtoutfile{MSC-class: \theprimaryclass\ifx\thesecondaryclass\relax\else, \thesecondaryclass\fi}
\immediate\write\gtoutfile{Journal-ref: Geom. Topol. \thevolumenumber
(\thevolumeyear) \startpage-\finishpage}
\immediate\write\gtoutfile{Comments: Published by Geometry and Topology at}
\immediate\write\gtoutfile{\s\s http://www.maths.warwick.ac.uk/gt/GTVol\thevolumenumber/paper\thepapernumber.abs.html}
\immediate\write\gtoutfile{\noexpand\\}
\immediate\write\gtoutfile{}
\ifx\theasciiabstract\relax
\immediate\write\gtoutfile{\theabstract}\else
\immediate\write\gtoutfile{\theasciiabstract}\fi
\immediate\write\gtoutfile{}
\immediate\write\gtoutfile{\noexpand\\}
\immediate\write\gtoutfile{}
\immediate\closeout\gtoutfile}}  

\def\maketitlepage{\maketitlep\makeheadfile}
\let\maketitle\maketitlepage


\def\ifplaintex{\expandafter\ifx\csname documentclass\endcsname\relax}


\ifplaintex 
\hoffset 14truemm
\voffset 31truemm
\else
\headsep 23pt
\footskip 35pt
\hoffset -4truemm
\voffset 12.5truemm
\fi

\expandafter\ifx\csname beginpicture\endcsname\relax
\expandafter\ifx\csname documentclass\endcsname\relax
\input pictex \else\font\fiverm=cmr5
\input prepictex \input pictex \input postpictex \fi\fi

\def\gt{{\mathsurround=0pt\it $\cal G\mskip-2mu$eometry \&\ 
$\cal T\!\!$opology}}        

\def\gtp{{\mathsurround=0pt\it $\cal G\mskip-2mu$eometry \&\ 
$\cal T\!\!$opology $\cal P\!$ublications}}  


\def\lognumber#1{\def\thelognumber{#1}}
\def\volumenumber#1{\def\thevolumenumber{#1}}
\def\papernumber#1{\def\thepapernumber{#1}}
\def\volumeyear#1{\def\thevolumeyear{#1}}

\def\pagenumbers#1#2{\def\startpage{#1}\def\finishpage{#2}}
\def\published#1{\def\publishdate{#1}}
\def\proposed#1{\def\theproposer{#1}}
\def\seconded#1{\def\theseconders{#1}}
\def\received#1{\def\receiveddate{#1}}

\def\accepted#1{\def\accepteddate{#1}}
\def\asciititle#1{\def\theasciititle{#1}}

\long\def\asciiabstract#1{\long\def\theasciiabstract{#1}}


\let\\\par\let\thelognumber\relax
\let\thevolumenumber\relax\let\thepapernumber\relax
\let\thevolumeyear\relax\let\thesamplenumber\relax\let\startpage\relax
\let\finishpage\relax\let\publishdate\relax\let\receiveddate\relax
\let\reviseddate\relax\let\accepteddate\relax\let\theasciititle\relax
\let\theasciiauthors\relax
\let\theasciiabstract\relax
\let\theasciiemail\relax\let\theshortauthors\relax\let\theshorttitle\relax

\long\def\maketitlep{   

\count0=\startpage

\gt\hfill      
\beginpicture
\setcoordinatesystem units <0.33truein, 0.33truein> point at 2.2 0.9
\setplotsymbol ({$\cal G$})
\plotsymbolspacing=9truept
\circulararc 315 degrees from 0 1 center at 0 0
\setplotsymbol ({$\cal T$})
\circulararc 315 degrees from 1 -1 center at 1 0
\endpicture
%
\break
{\small\ifx\thesamplenumber\relax 
Volume \else Sample
\fi\thevolumenumber\ (\thevolumeyear)
\startpage--\finishpage\nl
Published: \publishdate}
\vglue 0.5truein plus 0.4fil minus 0.1truein

{\parskip=0pt\leftskip 0pt plus 1fil\def\\{\par\smallskip}{\ifplaintex\large
\else\Large\fi\bf\thetitle}\par\medskip}   

\vglue 0pt plus 0.1fil 

{\parskip=0pt\leftskip 0pt plus 1fil\def\\{\par}{\sc\theauthors}
\par\medskip}

\vglue 0pt plus 0.1fil 

{\small\parskip=0pt\let\newline\\
{\leftskip 0pt plus 1fil\def\\{\par}{\sl\theaddress}\par}
\expandafter\ifx\theemail\relax    
\relax\else\vglue 5pt plus 0.02fil minus 2pt\def\\{\stdspace{\rm 
and}\stdspace} 
\cl{Email:\stdspace\tt\theemail}\fi
\ifx\theurl\relax                  
\relax\else\vglue 5pt plus 0.02fil minus 2pt\def\\{\stdspace{\rm 
and}\stdspace}
\cl{URL:\stdspace\tt\theurl}\fi\par}

\vglue 7pt plus 0.3fil minus 3pt

{\bf Abstract}
\vglue 5pt plus 0.1fil minus 2pt

\theabstract

\vglue 7pt plus 0.3fil minus 3pt

{\bf AMS Classification numbers}\quad Primary:\quad \theprimaryclass

Secondary:\quad \thesecondaryclass

\vglue 5pt plus 0.3fil minus 2pt

{\bf Keywords}\quad \thekeywords

\vglue 10pt plus 0.5fil minus 5pt

{\small  Proposed: \theproposer\hfill Received: \receiveddate\nl
Seconded: \theseconders\hfill 
\ifx\reviseddate\relax                         
Accepted: \accepteddate                        
\else
Revised: \reviseddate                          
\fi}
\eject
}       

\let\maketitlepage\maketitlep
\let\maketitle\maketitlepage


\font\phead=cmsl9 scaled 950
\font\lhead=cmsl9 scaled 1050
\font\pnum=cmbx10 scaled 913
\font\lnum=cmbx10 
\font\pfoot=cmsl9 scaled 950
\font\lfoot=cmsl9 scaled 1050
\ifplaintex
\headline{\vbox to 0pt{\vskip -4.5mm\line{\small\phead\ifnum
\count0=\startpage ISSN 1364-0380 (on line)
1465-3060 (printed) \hfill {\pnum\folio}\else\ifodd\count0\def\\{ }%
\ifx\theshorttitle\relax\thetitle\else\theshorttitle\fi\hfill{\pnum\folio}
\else\def\\{ and }{\pnum\folio}\hfill\ifx\theshortauthors\relax\theauthors
\else\theshortauthors\fi\fi\fi}\vss}}
\footline{\vbox to 0pt{\vglue 0mm\line{\small\pfoot\ifnum\count0=\startpage
\copyright\ \gtp\hfill\else
\gt, Volume \thevolumenumber\ (\thevolumeyear)\hfill\fi}\vss
}}
\else
\makeatletter
\def\@oddhead{{\small\lhead\ifnum\count0=\startpage ISSN 1364-0380 (on line)
1465-3060 (printed) \hfill {\lnum\number\count0}\else\ifodd\count0
\def\\{ }\ifx\theshorttitle\relax \thetitle \else\theshorttitle\fi\hfill
{\lnum\number\count0}\else\def\\{ and }{\lnum\number\count0}
\hfill\ifx\theshortauthors\relax 
\theauthors\else\theshortauthors\fi\fi\fi}}\def\@evenhead{\@oddhead}
\def\@oddfoot{\small\lfoot\ifnum\count0=\startpage\copyright\ \gtp\hfill\else
\gt, Volume \thevolumenumber\ (\thevolumeyear)\hfill\fi}
\def\@evenfoot{\@oddfoot}
\makeatother
\fi


\newwrite\gtoutfile
\long\gdef\makeheadfile{  
{\def\\{, }\def\s{ }
\immediate\openout\gtoutfile head.xxx
\immediate\write\gtoutfile{To: math@arxiv.org}
\immediate\write\gtoutfile{Subject: put or rep NNNNN:pppp}
\immediate\write\gtoutfile{--text follows this line--}
\immediate\write\gtoutfile{Proxy-for: \ifx\theasciiauthors\relax
\theauthors\else\theasciiauthors\fi\s<\ifx\theasciiemail\relax\theemail\else\theasciiemail\fi>}
\immediate\write\gtoutfile{\noexpand\\}
\immediate\write\gtoutfile{Authors: \ifx\theasciiauthors\relax
\theauthors\else\theasciiauthors\fi}
{\def\\{ }\immediate\write\gtoutfile{Title: \ifx\theasciititle\relax
\thetitle\else\theasciititle\fi}}
\immediate\write\gtoutfile{Subj-class: GT or SG or MG etc}
\immediate\write\gtoutfile{MSC-class: \theprimaryclass\ifx\thesecondaryclass\relax\else, \thesecondaryclass\fi}
\immediate\write\gtoutfile{Journal-ref: Geom. Topol. \thevolumenumber
(\thevolumeyear) \startpage-\finishpage}
\immediate\write\gtoutfile{Comments: Published by Geometry and Topology at}
\immediate\write\gtoutfile{\s\s http://www.maths.warwick.ac.uk/gt/GTVol\thevolumenumber/paper\thepapernumber.abs.html}
\immediate\write\gtoutfile{\noexpand\\}
\immediate\write\gtoutfile{}
\ifx\theasciiabstract\relax
\immediate\write\gtoutfile{\theabstract}\else
\immediate\write\gtoutfile{\theasciiabstract}\fi
\immediate\write\gtoutfile{}
\immediate\write\gtoutfile{\noexpand\\}
\immediate\write\gtoutfile{}
\immediate\closeout\gtoutfile}}  

\def\maketitlepage{\maketitlep\makeheadfile}
\let\maketitle\maketitlepage

\input epsf.tex
\useforwardrefs

\lognumber{250}
\volumenumber{6}
\papernumber{13}
\volumeyear{2002}
\pagenumbers{393}{401}
\received{30 April 2001}
\accepted{9 July 2002}
\published{21 July 2002}
\proposed{Robion Kirby}
\seconded{Wolfgang Metzler, Ronald Stern}

\reflist

\key{A20}	{\bf J\,W Alexander}, {\it Note on Riemann spaces},
Bull. Amer. Math. Soc. 26 (1920) 370--373

\key{BE78}	{\bf I Berstein}, {\bf A\,L Edmonds}, {\it The degree
and branch set of a branched covering}, Invent. Math. 45 (1978)
213--220

\key{F57}	{\bf R\,H Fox}, {\it Covering spaces with
singularities}, Algebraic Geometry and Topology, A symposium in
honour of S. Lefschetz, Princeton 1957, 243--257

\key{Hi74}	{\bf H\,M Hilden}, {\it Every closed orientable
3--manifold is a 3--fold branched covering space of $S^3$}, 
Bull. Amer. Math. Soc. 80 (1974) 1243--1244

\key{HM80}	{\bf H\,M Hilden}, {\bf J\,M Montesinos}, 
{\it Lifting surgeries to branched covering spaces}, 
Trans. Amer. Math. Soc. 259 (1980) 157--165

\key{Hr74}	{\bf U Hirsch}, {\it \"Uber offene Abbildungen auf die
3--Sph\"are}, Math. Z. 140 (1974) 203--230

\key{HK79}	{\bf F Hosokawa}, {\bf A Kawauchi}, {\it Proposals
for unknotted surfaces in four-spaces}, Osaka J. Math. 16
(1979) 233--248

\key{K89}	{\bf S Kamada}, {\it Non-orientable surfaces in 4--space}, 
Osaka J. Math. 26 (1989) 367--385

\key{K89}	{\bf R Kirby}, {\it Problems in low-dimensional
topology}, Stud. Adv. Math. 2, part 2, Amer. Math. Soc.
and International Press 1997

\key{M69}	{\bf W\,S Massey}, {\it Proof of a conjecture of
Whitney}, Pacific J. Math. 31 (1969) 143--156

\key{M74}	{\bf J\,M Montesinos}, {\it A representation of closed,
orientable 3--manifolds as 3--fold branched coverings of $S^3$},
Bull. Amer. Math. Soc. 80 (1974) 845--846

\key{M78}	{\bf J\,M Montesinos}, {\it 4--manifolds, 3--fold covering
spaces and ribbons}, Trans.\break Amer. Math. Soc. 245 (1978) 453--467

\key{P95} {\bf R Piergallini}, {\it Four-manifolds as $4$--fold 
branched covers of $S^4$}, Topology 34 (1995) 497--508

\key{V84} {\bf O\,Ja Viro}, {\it Signature of branched covering},
Trans. Mat. Zametki 36 (1984) 549--557

\endreflist

\title{4--manifolds as covers of the 4--sphere\\%
       branched over non-singular surfaces}                   
\authors{Massimiliano Iori \\ Riccardo Piergallini}                  
\address{Dipartimento di Matematica e Informatica \\
         Universit\`a di Camerino -- Italia}                  

\asciititle{4-manifolds as covers of the 4-sphere 
branched over non-singular surfaces}                   

\email{riccardo.piergallini@unicam.it}                     

\abstract
We prove the long-standing Montesinos conjecture that any closed
oriented PL 4--manifold $M$ is a simple covering of $S^4$
branched over a locally flat surface (cf~[\ref{M78}]).
In fact, we show how to eliminate all the node singularities of
the branching set of any simple 4--fold branched covering $M \to
S^4$ arising from the representation theorem given in
[\ref{P95}]. Namely, we construct a suitable cobordism between the
5--fold stabilization of such a covering (obtained by adding a
fifth trivial sheet) and a new 5--fold covering $M \to S^4$
whose branching set is locally flat. It is still an open
question whether the fifth sheet is really needed or not.
\endabstract

\asciiabstract{
We prove the long-standing Montesinos conjecture that any closed
oriented PL 4-manifold M is a simple covering of S^4 branched over a
locally flat surface (cf [J M Montesinos, 4-manifolds, 3-fold
covering spaces and ribbons, Trans. Amer. Math. Soc. 245 (1978)
453--467]).  In fact, we show how to eliminate all the node
singularities of the branching set of any simple 4-fold branched
covering M \to S^4 arising from the representation theorem given in [R
Piergallini, Four-manifolds as 4-fold branched covers of S^4, Topology
34 (1995) 497--508]. Namely, we construct a suitable cobordism between
the 5-fold stabilization of such a covering (obtained by adding a
fifth trivial sheet) and a new 5-fold covering M --> S^4 whose
branching set is locally flat. It is still an open question whether
the fifth sheet is really needed or not.}

\primaryclass{57M12}
\secondaryclass{57N13}             
\keywords{4--manifolds, branched coverings, locally flat branching surfaces}
\maketitlepage

\section{Introduction}

The idea of representing manifolds as branched covers of
spheres, extending the classical theory of ramified surfaces
introduced by Riemann, is due to Alexander [\ref{A20}] and dates
back to 1920. He proved that for any orientable closed PL
manifold $M$ of dimension $m$ there is a branched covering of
$M \to S^m$.

We recall that a non-degenerate PL map $p\co M \to N$ between
compact PL manifolds is called a {\it branched covering} if there
exists an $(m-2)$--subcomplex $B_p \subset N$, the {\it branching
set} of $p$, such that the restriction $p_|\co  M - p^{-1}(B_p) \to
N - B_p$ is an ordinary covering of finite degree $d$.
If $B_p$ is minimal with respect to such property, then we have
$B_p = p(S_p)$, where $S_p$ is the {\it singular set} of
$p$, that is the set of points at which $p$ is not locally
injective. In this case, both $B_p$ and $S_p$, as well as the 
{\it pseudo-singular set} $S'_p = {\rm Cl}(p^{-1}(B_p)-S_p)$,
are (possibly empty) homogeneously $(m-2)$--dimensional
complexes.

Since $p$ is completely determined (up to PL homeomorphism) by
the ordinary covering $p_|$ (cf~[\ref{F57}]), we can describe it
in terms of its branching set $B_p$ and its monodromy
$\omega_p \co  \pi_1(N-B_p) \to \Sigma_{d}$ (uniquely defined up
to conjugation in $\Sigma_{d}$, depending on the numbering of
the sheets). 

If $N = S^m$ then a convenient description of $p$ can be given
by labelling each $(m-2)$--simplex of $B_p$ by the monodromy of
the corresponding meridian loop, since such loops generate the
fundamental group $\pi_1(S^m-B_p)$.

Therefore, we can reformulate the Alexander's result as follows:
any orientable closed PL manifold $M$ of dimension $m$ can be
represented by a labelled $(m-2)$--subcomplex of $S^m$.

Of course, in order to make such representation method
effective, some control is needed on the degree $d$ and on the
complexity of the local structure of $B_p$ and $\omega_p$.
Unfortunately, there is no such control in the original
Alexander's proof, being $d$ dependent on the number of
simplices of a triangulation of $M$ and $B_p$ equal to the
$(m-2)$--skeleton of an $m$--simplex. Even at the present, as
far as we know, the only general (for any $m$) results in
this direction are the negative ones obtained by Berstein and
Edmonds [\ref{BE78}]: for representing all the $m$--manifolds at
least $m$ sheets are necessary (for example this happens of the
$m$--torus $T^m$) and in general we cannot require $B_p$ to be
non-singular (the counterexamples  they give have dimension $m 
\geq 8$). On the contrary, the situation is much better for $m 
\leq 4$.

The case of surfaces is trivial: the closed (connected)
orientable surface $T_g$ of genus $g$ is a 2--fold cover of
$S^2$ branched over $2g + 2$ points. For $m = 3$, Hilden
[\ref{Hi74}], Hirsch [\ref{Hr74}] and Montesinos [\ref{M74}]
independently proved that any orientable closed (connected)
3--manifold is a simple 3--fold cover of $S^3$ branched over a
knot.

For $m = 4$, the representation theorem proved by Piergallini
[\ref{P95}] asserts that any orientable closed (connected) PL
$4$--manifold is a simple $4$--fold cover of $S^4$ branched over a
transversally immersed PL surface. Simple means that the
monodromy of each meridian loop is a transposition. On the other
hand, a transversally immersed PL surface is a subcomplex which
is a locally flat PL surface at all its points but a finite
number of {\it nodes} (transversal double points).
So, the local models (up to PL equivalence) for the labelled
branching set are the ones depicted in Figure
\ref{localmodels/fig}, where $\{i,j,k,l\} = \{1,2,3,4\}$ (the
monodromies of the meridian loops corresponding to sheets of the
branching set meeting at a node must be disjoint). We remark
that in general the branching surface cannot be required to be
orientable (cf~[\ref{P95}], [\ref{V84}]).

\figure
\epsfbox{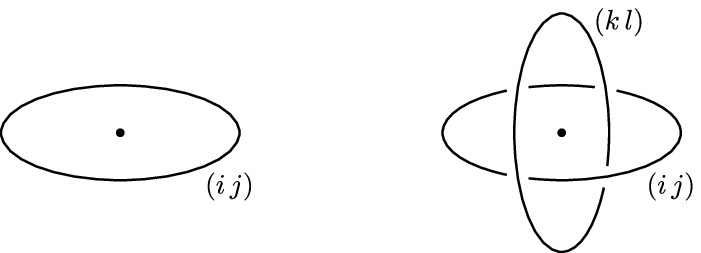}
\key{localmodels/fig}
\endfigure

The question whether the nodes can be eliminated in order to get
non-singular branching surfaces, as proposed by Montesinos in
[\ref{M78}], was left open in [\ref{P95}]. 

In the next section we show how elimination of nodes can be
performed up to cobordism of coverings, after the original
4--fold covering has been stabilized by adding a fifth trivial
sheet. This proves the following representation theorem.

\proclaim{Theorem}
Any orientable closed (connected) PL $4$--manifold is a simple
$5$--fold cover of $S^4$ branched over a locally flat PL surface.
\endproc

\section{Elimination of nodes}
\key{nodes/sec}

Let $M$ be an orientable closed (connected) PL 4--manifold and
let $p\co M \to S^4$ be a 4--fold covering branched over a
transversally immersed PL surface $F \subset S^4$ given by
Theorem B of [\ref{P95}]. We denote by $q\co  M \to S^4$ the 5--fold
branched covering obtained by stabilizing $p$ with an extra
trivial sheet. In terms of labelled branching set this means
adding to the surface $F$, labelled with transpositions in
$\Sigma_4$, a separate unknotted 2--sphere $S$ labelled with the
transposition $(4\,5)$, as schematically shown in Figure
\ref{stabilization/fig}.

\figure
\epsfbox{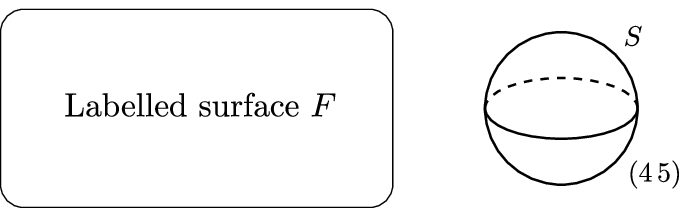}
\key{stabilization/fig}
\endfigure

Looking at the proof of Theorem B of [\ref{P95}], we see that
nodes of the branching set of $p$ come in pairs, in such a way
that each pair consists of the end points of a simple arc
contained in $F$ and all these arcs are disjoint from each
other.

Let $\alpha_1, \dots, \alpha_n \subset F$ be such arcs
and let $\nu_i$ and $\nu'_i$ be the nodes joined by $\alpha_i$. 
The intersection of $F \cup S$  with a sufficiently small regular 
neighborhood $N(\alpha_i)$ of $\alpha_i$ in $S^4$ consists of a 
disk $A_i$ containing $\alpha_i$ and two other disks $B_i$ and 
$B'_i$ transversally meeting $A_i$ respectively at $\nu_i$ and 
$\nu'_i$. Up to labelled isotopy, we can assume $B_i$ and $B'_i$ 
labelled with $(1\,2)$ and $A_i$ labelled with $(3\,4)$, as in 
Figure \ref{alpha/fig} (remember that the monodromy of $p$ is
transitive, since $M$ is connected). We also assume the
$N(\alpha_i)$'s disjoint from each other. 

\figure
\epsfbox{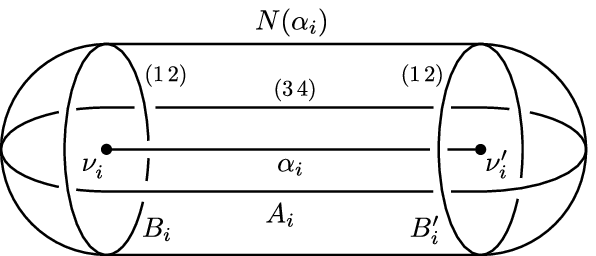}
\key{alpha/fig}
\endfigure

For future use, we modify the branching surface $F \cup S$ by
``finger move'' labelled isotopies, in order to introduce inside
each $N(\alpha_i)$ two more small trivial disks $C_i$ and $C'_i$
respectively labelled by $(2\,4)$ and $(4\,5)$, as shown in
Figure \ref{fingers/fig}. This modification has the effect of
connecting $q^{-1}(N(\alpha_i))$ making it PL equivalent to $S^1
\times B^3$.

\figure
\epsfbox{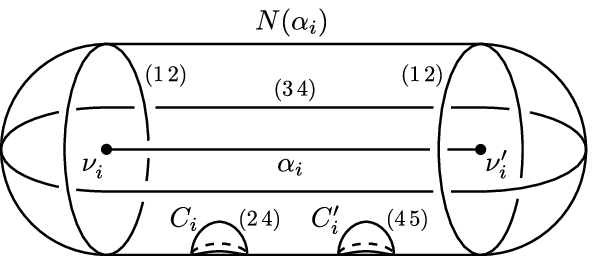}
\key{fingers/fig}
\endfigure

Now, we consider the orientable 5--manifold $T = S^4 \times [0,1]
\cup H_1 \cup \dots \cup H_n$ obtained by attaching to $S^4
\times [0,1]$ a 1--handle $H_i$ for each pair of nodes $\nu_i$,
$\nu'_i$. The attaching cells of each $H_i$ are $N(\nu_i) \times
\{1\}$ and $N(\nu'_i) \times \{1\}$, where $N(\nu_i)$ and
$N(\nu'_i)$ are regular neighborhoods $\nu_i$ and $\nu'_i$ in
$N(\alpha_i) - (C_i \cup C'_i)$, such that all the intersections
$D_i = N(\nu_i) \cap A_i$, $E_i = N(\nu_i) \cap B_i$, $D'_i = 
N(\nu'_i) \cap A_i$ and $E'_i = N(\nu'_i) \cap B'_i$ are again 
disks.

The product covering $q \times {\rm id}_{[0,1]}\co  M \times [0,1]
\to S^4 \times [0,1]$ can be extended to a new 5--fold simple
branched covering $r\co  W \to T$, where $W$ is the result of
adding appropriate $1$--handles to $M\times [0,1]$ over the
$H_i's$. In fact, the restrictions of $q \times \{1\}$ over 
$N(\nu_i) \times \{1\}$ and $N(\nu'_i) \times \{1\}$ are 
equivalent, hence, by a suitable choice of the attaching map of 
$H_i$, we can define $r$ over $H_i \cong B^4 \times [0,1]$ just 
by crossing the first restriction with the identity of $[0,1]$. 
Namely, the pair $(H_i,B_r \cap H_i)$ is equivalent to 
$(N(\nu_i),D_i \cup E_i) \times [0,1]$, with the monodromy of 
the meridian loops around $D_i \times [0,1]$ and 
$E_i \times [0,1]$ respectively equal to $(1\,2)$ and $(3\,4)$.
Then, $r^{-1}(H_i)$ consists of three 1--handles attached to $M
\times [0,1]$ at the three pairs of 4--cells making up the pair
$(q^{-1}(N(\nu_i)),q^{-1}(N(\nu'_i))) \times \{1\}$. 
We denote by $H'_i$, $H''_i$ and $H'''_i$ these 1--handles in
such a way that they involve respectively the sheets 1 and 2,
the sheets 3 and 4, and the sheet 5 (see Figure
\ref{handles1/fig}, where the lighter lines represent the
pseudo-singular set).  We remark that the branching set $B_r$ is
a locally flat PL 3--manifold at all points but one transversal
double arc inside each $H_i$ between $\nu_i \times \{1\}$ and
$\nu'_i \times \{1\}$.

\figure
\epsfbox{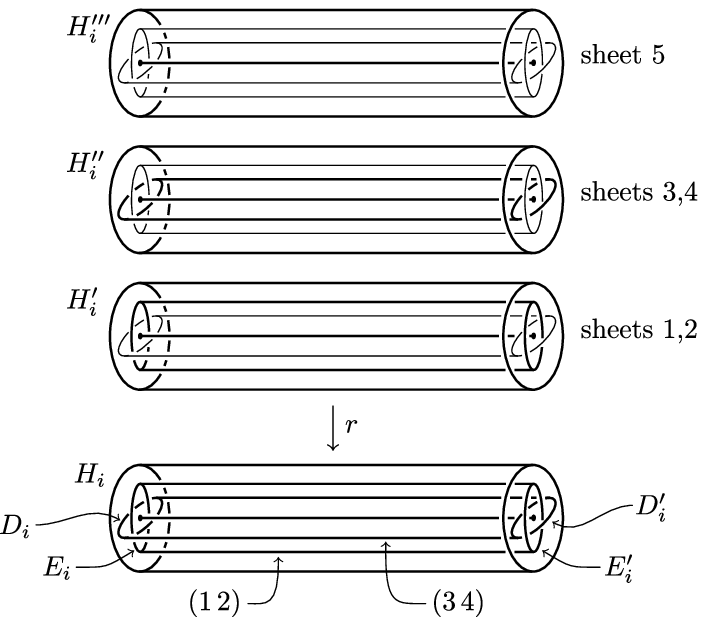}
\key{handles1/fig}
\endfigure

At this point, we want to simultaneously attach to $T$ and $W$
some 2--handles in order to kill the 1--handles $H_1, \dots, H_n$
attached to $S^4 \times [0,1]$ and the 1--handles $H'_1, H''_1,
H'''_1, \dots, H'_n, H''_n, H'''_n$ attached to $M \times
[0,1]$, taking care that the branched covering $r$ can be
extended to these 2--handles.

For each $i = 1, \dots, n$, we consider a simple loop $\lambda_i$ 
inside ${\rm Bd} T \cap (N(\alpha_i) \times \{1\} \cup H_i) - B_r$
running through $H_i$ once and linking both the disks $C_i \times 
\{1\}$ and $C'_i \times \{1\}$ once, as shown in Figure 
\ref{lambda/fig}. We observe that $r^{-1}(\lambda_i)$ consists of 
three loops $\lambda'_i,\lambda''_i,\lambda'''_i \subset {\rm Bd}
W - (S_r \cup S'_r)$, such that: $\lambda'_i$ runs through 
$H'_i$ once and avoids $H''_i \cup H'''_i$, $\lambda''_i$
runs through $H''_i$ once and avoids $H'_i \cup H'''_i$, while
$\lambda'''_i$ runs through each of $H'_i$, $H''_i$ and $H'''_i$
once. 

\figure
\epsfbox{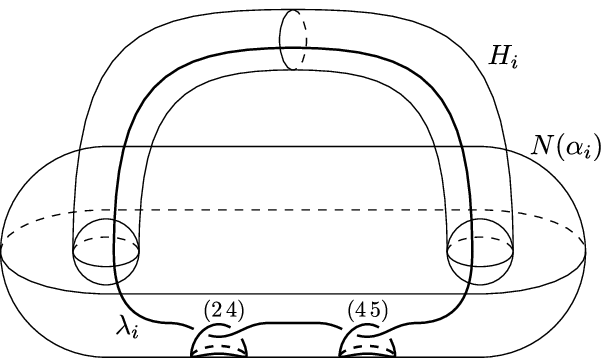}
\key{lambda/fig}
\endfigure

Then, the 5--manifold $T \cup L_1 \cup \dots \cup L_n$ obtained
by attaching to $T$ the 2--handle $L_i$ along each loop
$\lambda_i$ (with arbitrary framing), is PL homeomorphic to $S^4
\times [0,1]$, since each $L_i$ kills the corresponding $H_i$.

Analogously, the 5--manifold $W \cup (L'_1 \cup L''_1 \cup
L'''_1) \cup \dots \cup (L'_n \cup L''_n \cup L'''_n)$
obtained by attaching to $W$ the 2--handles $L'_i$, $L''_i$
and $L'''_i$ along the loops $\lambda'_i$, $\lambda''_i$ and
$\lambda'''_i$ (with arbitrary framings), is PL homeomorphic to
$M \times [0,1]$. In fact, we can cancel first each $L'''_i$
with the corresponding $H'''_i$ and then each $L'_i$ and $L''_i$
respectively with $H'_i$ and $H''_i$.

By choosing the attaching framings of the 2--handles $L'_i$,
$L''_i$ and $L'''_i$ accordingly with the ones of the 2--handle
$L_i$, we can extend the covering $r$ to such 2--handles as
suggested by Figure \ref{handles2/fig}, where the branching set
consists of the labelled 3--cells $F_i$ and $G_i$ transversal to
the 2--handle $L_i$. Namely, we can glue the covering represented
in the figure with $r$, since they coincide over the attaching 
tube around $\lambda_i$. Then, we can identify $L_i'$ and $L_i''$ 
respectively with the trivial components over $L_i$ corresponding 
to sheets 1 and 3, and $L_i'''$ with the non-trivial component 
over $L_i$ corresponding to sheets $2$, $4$ and $5$.

\figure
\epsfbox{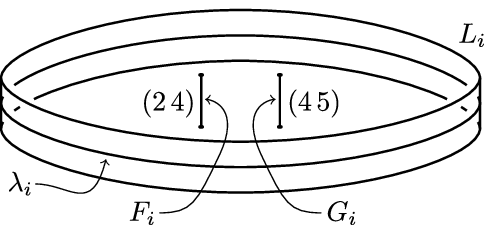}
\key{handles2/fig}
\endfigure

In this way, we get an extension of $r$ which is PL equivalent to
a new branched covering $s\co  M \times [0,1] \to S^4 \times
[0,1]$. Up to the natural identification between fibers and
factors, the restriction of $s$ over $S^4 \times \{0\}$
coincides with $q$, while the restriction over $S^4 \times
\{1\}$ gives us a new 5--fold simple branched covering $q'\co M \to
S^4$.

The branching set $B_{q'}$ of $q'$ is a locally flat PL surface
in $S^4$. In fact, it is isotopically equivalent to the result of
the following modifications performed on $B_q = F \cup S$, due
to attaching handles: for each $i = 1, \dots, n$, the disks
$D_i$, $D'_i$, $E_i$ and $E'_i$ are replaced by linked pipes
respectively connecting ${\rm Bd} D_i$ with ${\rm Bd} D'_i$ and 
${\rm Bd} E_i$ with ${\rm Bd} E'_i$; for each $i = 1, \dots, n$, 
the new trivial spheres ${\rm Bd} F_i$ and ${\rm Bd} G_i$ are 
added on.

\section{Final remarks}
\key{remark/sec}

The argument used in the previous section for eliminating nodes,
with some minor variation, allows us to perform a variety of
different modifications on branched coverings.

We can eliminate any pair of isolated singularities of the
branching set, which are equivalent up to orientation reversing
PL homeomorphisms, provided that the covering has at least one
sheet more than the ones involved in them. For instance, this
is a way, alternative with respect to the one of [\ref{P95}], 
to remove cusps from the branching set of a simple 4--fold 
covering of $S^4$.

On the other hand, by choosing the attaching balls of the
1--handle $H_i$ centred at two non-singular points of the
branching set with the same monodromy and letting the
attaching loop of the 2--handle $L_i$ have trivial monodromy, we
get a new approach to surgery of simple branched coverings
along symmetric knots (see [\ref{M78}]). In fact, in
this case we have $d-1$ handles over $H_i$ and $d$ handles over
$L_i$, where $d$ is the degree of the covering, and
after cancellation we are left with one 2--handle attached
to the covering manifold along the unique loop in the
counterimage of the arc $\alpha_i$. Surgeries of greater indices
(see [\ref{HM80}]) can be realized similarly.

With a different choice of the monodromies, we can also perfom
surgeries on the branching set without changing the covering
manifold up to PL homeomorphisms. In particular, we get the move
shown in Figure \ref{move/fig}, which is the double of the move
in Figure 12 of [\ref{P95}].

\figure
\epsfbox{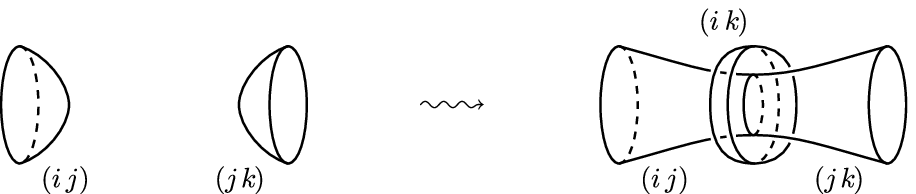}
\key{move/fig}
\endfigure

By using this move, we can connect all the non-trivial
components of the branching surface, provided that the degree
of the covering is at least $3$, in such a way that the
branching surface of the theorem can be assumed to have the 
following special form: $F = G \cup S_1 \cup \dots \cup S_k$, 
where $G \subset S^4$ is connected and $S_1, \dots, S_k$ is a 
family of separate trivial 2--spheres. Furthermore, we can perform
hyperbolic transformations of $G$ in order to make it unknotted
(cf~[\ref{HK79}], [\ref{K89}]).

We observe that, in some sense, $G$ represents the cobordism class
of the covering manifold $M$, being $\sigma(M) = - F \cdot F/2 =
-G\cdot G/2$ (cf~[\ref{V84}]). On the other hand, the $S_i$'s
cannot be eliminated in general, that is the branching surface
cannot be required to be connected. In fact, given any covering 
$M \to S^4$ branched over a locally flat PL surface $F$, we
have $\chi(M) = 2d - \chi(F)$, where $d$ is the degree of the
covering. Then, by the Whitney inequality for the
self-intersection of non-orientable surfaces in $S^4$ 
(cf~[\ref{M69}]), $F$ must have at least $d + |\sigma(M)|/2
- \chi(M)/2$ components.

\medskip

Finally, we remark that our argument heavily depends on the
fifth extra sheet for the elimination of nodes, hence it seems
useless for solving the following question that remains still
open (cf~Problem 4.113 of Kirby's problem list [\ref{K89}]):

\proclaim{Question}
Is any orientable closed (connected) PL $4$--manifold a simple
$4$--fold cover of $S^4$ branched over a locally flat PL surface?
\endproclaim

\vfill\break

\references 

\bye